\documentclass{amsart}% option [reqno] puts equat numb to the right

\usepackage[mathscr]{eucal}% So-called Euler fonts, allows \mathscr
\usepackage{amssymb}% allows special arrows like >-> e.g.
\usepackage[usenames,dvipsnames]{color}% allows colors
\usepackage{amsthm}% allows theorem* , e.g.
\usepackage{bbold}% allows \unit \mathbb{1}
\usepackage{enumerate}% allows easy change of label
\usepackage{array}% allows array
\usepackage{amsmath}% allows pmatrix

\usepackage[colorlinks=true,linkcolor={Brown},citecolor={Brown},urlcolor={Brown}]{hyperref}

\usepackage[all]{xy}
\SelectTips{cm}{}
\newdir{ >}{{}*!/-10pt/\dir{>}}

\numberwithin{equation}{section}% makes equat numb contain the section
\swapnumbers %pour que le numéro apparaisse devant le théorème

\newtheorem{Thm}[equation]{Theorem}
\newtheorem*{Thm*}{Theorem}
\newtheorem{Cor}[equation]{Corollary}
\theoremstyle{remark}
\newtheorem{Def}[equation]{Definition}
\newtheorem{Not}[equation]{Notation}
\newtheorem{Exa}[equation]{Example}
\newtheorem{Assum}[equation]{Assumption}

\newtheorem{Rem}[equation]{Remark}
\newtheorem{Point}[equation]{Point}
\newcommand{\bbA}{\mathbb{A}}
\newcommand{\bbC}{\mathbb{C}}
\newcommand{\bbF}{\mathbb{F}}

\newcommand{\bbP}{\mathbb{P}}
\newcommand{\bbQ}{\mathbb{Q}}

\newcommand{\bbZ}{\mathbb{Z}}
\newcommand{\cat}[1]{\mathscr{#1}}%or: \newcommand{\cat}[1]{\mathcal{#1}}
\newcommand{\eg}{{\sl e.g.}}
\newcommand{\Endcat}[1]{\End_{\cat #1}}

\newcommand{\ie}{{\sl i.e.}\ }
\newcommand{\isoto}{\overset{\sim}{\,\to\,}}

\newcommand{\SET}[2]{\big\{\,#1\,\big|\,#2\,\big\}}
\newcommand{\too}{\mathop{\longrightarrow}\limits}
\newcommand{\unit}{\mathbb{1}}% unit for \otimes
\newcommand{\potimes}[1]{^{\otimes #1}}% tensor power
\newcommand{\aka}{{a.\,k.\,a.}\ }
\newcommand{\SHc}{\SH^c}
\newcommand{\SHcp}{\SH^c_{(p)}}
\newcommand{\Db}{\Der^{\smallb}}% derived bounded category
\newcommand{\adh}[1]{\overline{#1}}% adherence
\newcommand{\adhpt}[1]{\adh{\{#1\}}}% adherence of a pt
\newcommand{\Dperf}{\Der^{\smallperf}}% derived category of perfect compl
\newcommand{\Kb}{\Komp^{\smallb}}% htpy bounded category
\newcommand{\grmod}{\,\text{--}\grmodname}
\newcommand{\pproj}{\,\text{--}\proj}% category of projective modules
\newcommand{\stmod}{\,\text{--}\stmodname}%
\newcommand{\strel}{\,\text{--}\strelname}%
\newcommand{\mmod}{\,\text{--}\modname}%
\newcommand{\sbull}{{\scriptscriptstyle\bullet}}
\newcommand{\Spech}{\Spec^{\text{\rm h}}}
\newcommand{\ideal}[1]{\langle #1\rangle}
\newcommand{\SpcK}{\Spc(\cat K)}% most used
\newcommand{\gp}{\mathfrak{p}}% prime p
\newcommand{\gq}{\mathfrak{q}}% prime q
\newcommand{\Aone}{{\bbA^{\!\scriptscriptstyle 1}}}

\newcommand{\Aut}{\operatorname{Aut}}
\newcommand{\coh}{\operatorname{coh}}
\newcommand{\cone}{\operatorname{cone}}

\newcommand{\Der}{\operatorname{D}}
\newcommand{\DM}{\operatorname{DM}}
\newcommand{\DTM}{\operatorname{DTM}}
\newcommand{\End}{\operatorname{End}}
\newcommand{\et}{\operatorname{et}}
\newcommand{\grmodname}{\operatorname{grmod}}
\newcommand{\Hom}{\operatorname{Hom}}
\newcommand{\Ind}{\operatorname{Ind}}
\newcommand{\Ker}{\operatorname{Ker}}
\newcommand{\Komp}{\operatorname{K}}
\newcommand{\modname}{\operatorname{mod}}
\newcommand{\mot}{\operatorname{mot}}
\newcommand{\proj}{\operatorname{proj}}
\newcommand{\Proj}{\operatorname{Proj}}
\newcommand{\Res}{\operatorname{Res}}
\newcommand{\rmH}{\operatorname{H}}
\newcommand{\rmL}{\operatorname{L}}
\newcommand{\rk}{\operatorname{rk}}
\newcommand{\SH}{\operatorname{SH}}
\newcommand{\smallb}{\operatorname{b}}
\newcommand{\smallperf}{\operatorname{perf}}
\newcommand{\Spc}{\operatorname{Spc}}
\newcommand{\Spec}{\operatorname{Spec}}
\newcommand{\stmodname}{\operatorname{stmod}}
\newcommand{\strelname}{\operatorname{strel}}
\newcommand{\supp}{\operatorname{supp}}

\begin{document}
%----------------------------------------------------------------------

\title{A guide to tensor-triangular classification}
\author{Paul Balmer}
\date{2019 December 17}

\address{Paul Balmer, Mathematics Department, UCLA, Los Angeles, CA 90095-1555, USA}
\email{balmer@math.ucla.edu}
\urladdr{http://www.math.ucla.edu/$\sim$balmer}

\begin{abstract}
This is a chapter of the Handbook of Homotopy Theory, that surveys the classifications of thick tensor-ideals.
\end{abstract}


\thanks{Research supported by NSF grant~DMS-1600032.}

\maketitle

%------------------------------------------------------------------------------
%----------------------------------------------------------------------
\section{Introduction}
\label{sec:intro}%
%----------------------------------------------------------------------

Stable homotopy theory shines across pure mathematics, from topology to analysis, from algebra to geometry. While its liturgy invokes Quillen model structures and $\infty$-categories, profane users around the world often speak the vernacular of \emph{triangulated categories}, as we shall do in this chapter.

Perhaps the first salient fact about stable homotopy categories is that in almost all cases they turn out to be \emph{wild} categories -- beyond the trivial examples of course. Dade famously began his paper~\cite{Dade78a} with the admonition ``There are just too many modules over p-groups!" and this truth resonates in all other fields as well: no hope to classify topological spaces up to stable homotopy equivalence; no more hope with complexes of sheaves, nor with equivariant $C^*$-algebras, nor with motives, etc, etc. One might dream that things improve with `small' objects (compact, rigid, or else) but the problem persists even there: Stable homotopy theory is just too complicated!

Faced with the complexity of stable homotopy categories, we are led to the following paradigm shift. A classification \emph{up to isomorphism} makes sense in any category, \ie as soon as we can speak of \emph{isomorphism}. But stable homotopy categories are more than mere categories: They carry additional structures, starting with the \emph{triangulation}. In the case of  a \emph{tensor}-triangulated category (\emph{tt-category} for short), as we consider in this chapter, we have two basic tools at hand: triangles and tensor. Instead of ignoring these additional structures, we should include them in the concept of
\begin{center}
\emph{tt-classification}
\index{classification of thick tensor-triangulated ideals}%
\end{center}
which is our nickname for \emph{classification up to the tensor-triangular structure}.

More precisely, we want to decide when two objects $X$ and $Y$ can be obtained from one another by using tensor with anything, direct sums, summands, cones, suspension, etc. In mathematical terms, we ask when $X$ and $Y$ generate the same thick triangulated tensor-ideals. Heuristically, if you can build $Y$ out of~$X$ by using the tt-structure then $X$ contains at least as much information as~$Y$. If you can go back and forth between~$X$ and~$Y$, then they contain the same amount of information.

The remarkable gain is that the tt-classification of an essentially small (rigid) tensor-triangulated category can \emph{always} be achieved by means of a geometric object, more precisely a spectral topological space, called its \emph{tensor-triangular spectrum}.\,(\footnote{\,Our use of the word `spectrum' comes from commutative algebra, as in the `Zariski spectrum', and should not be confused with the suspension-inverting `spectra' of topology.}) Let us highlight this starting point:

\smallbreak
\noindent\textbf{Fundamental fact:} \emph{Although almost every symmetric monoidal stable homotopy category~$\cat{K}$ is `wild' as a category, we always have a tt-classification of its objects, via a topological space, $\SpcK$, called the spectrum of~$\cat{K}$}.\,(\footnote{\,$\SpcK$ is a space in the universe containing the `set' of isomorphism classes of~$\cat{K}$.})
\smallbreak

This chapter is dedicated to a survey of tt-classifications across different examples, as far as they are known to the author at this point in time.

\medbreak

The original idea of classifying objects up to the ambient structure was born in topology, around Ravenel's conjectures~\cite{Ravenel84} and the `chromatic' theorems of Devinatz-Hopkins-Smith~\cite{DevinatzHopkinsSmith88,HopkinsSmith98}; this relied on Morava's work, among many other contributions. The ground-breaking insight of transposing from topology to other fields began with Hopkins~\cite{Hopkins87}. It is arguably Thomason~\cite{Thomason97} who first understood how essential the tensor was in the global story. We recall in Remark~\ref{rem:no-tensor-no-cry} why such a \emph{geometric} classification does not exist for mere triangulated categories, \ie without the tensor.

The tt-spectrum was introduced in~\cite{Balmer05a} and is reviewed in Section~\ref{sec:spc-classif}. The survey begins in Section~\ref{sec:topology}, with the initial example of topological stable homotopy theory. Section~\ref{sec:AG} touches commutative algebra and algebraic geometry. Section~\ref{sec:modular} is dedicated to stable module categories in modular representation theory and beyond. Section~\ref{sec:equivariant} discusses equivariant stable homotopy theory and Kasparov's equivariant $KK$-theory. Section~\ref{sec:motives} pertains to motives and $\Aone$-homotopy theory.

Everywhere, we have tried to give some idea of the actual tt-categories which come into play. When the amount of specialized definitions appears too high for this chapter, we simply point to the bibliographical references.

Finally, let us say a word about the bigger picture. In commutative algebra, the Zariski spectrum is not meant to be explicitly computed for every single commutative ring in the universe; instead, it serves as a stepping stone towards the geometric reasonings of algebraic geometry. In the same spirit, the tt-spectrum opens up a world of mathematical investigation, called \emph{tensor-triangular geometry}, which reaches far beyond classical algebraic geometry into the broad kingdom of stable homotopy theory. The short final Section~\ref{sec:tt-geometry} points to further reading in that direction.

\medbreak
\noindent\textbf{Acknowledgements}: I would like to thank Tobias Barthel, Ivo Dell'Ambrogio, Martin Gallauer, Beren Sanders and Greg Stevenson for their help in assembling this survey. I apologize to anyone whose tt-geometric results are not mentioned here: For the sake of pithiness, I chose to restrict myself to the topic of tt-classification.

\goodbreak

%----------------------------------------------------------------------
\section{The tt-spectrum and the classification of tt-ideals}
\label{sec:spc-classif}%
%----------------------------------------------------------------------

%
\begin{Def}
\index{tensor-triangulated category}%
\index{tt-category}%
A \emph{tt-category}, short for \emph{tensor-triangulated category}, is a triangulated category~$\cat{K}$ together with a symmetric monoidal structure
\[
\otimes\colon \cat{K}\times\cat{K}\too \cat{K}
\]
which is exact in each variable. See details in~\cite[App.\,A]{HoveyPalmieriStrickland97} or~\cite{Neeman01,KellerNeeman02}. The $\otimes$-unit is denoted~$\unit$.
\end{Def}

\begin{Assum}
\label{ass:small}%
Unless otherwise stated, we always assume that $\cat{K}$ is \emph{essentially small}, \ie has a set of isomorphism classes of objects. Subcategories $\cat{J}\subseteq\cat{K}$ are always assumed full and replete (\ie closed under isomorphisms).
\end{Assum}

\begin{Def}
\label{def:general-tt}%
A \emph{triangulated} subcategory $\cat{J}\subseteq \cat{K}$ is a non-empty subcategory such that whenever $X\to Y\to Z\to \Sigma X$ is an exact triangle in~$\cat{K}$ and two out of $X$, $Y$ and~$Z$ belong to~$\cat{J}$ then so does the third. A \emph{thick} subcategory $\cat{J}\subseteq\cat{K}$ is a triangulated subcategory closed under direct summands\,: if  $X\oplus Y\in\cat{J}$ then $X,Y\in\cat{J}$. A \emph{tt-ideal} $\cat{J}\subseteq\cat{K}$, short for \emph{thick tensor-ideal}, is a thick subcategory closed under tensoring with any object\,: $\cat{K}\otimes\cat{J}\subseteq\cat{J}$. A tt-ideal~$\cat{J}\subseteq\cat{K}$ is called \emph{radical} if $X\potimes{n}\in\cat{J}$ for~$n\ge 2$ forces $X\in\cat{J}$.
\end{Def}

\begin{Rem}
\label{rem:rigid}%
When every object of~$\cat{K}$ is rigid (\ie admits a dual, \aka \emph{strongly dualizable} \cite[\S\,2.1]{HoveyPalmieriStrickland97}), then we say that $\cat{K}$ is \emph{rigid} and we can show that every tt-ideal~$\cat{J}$ is automatically radical. See~\cite[Prop.\,2.4]{Balmer07}. So, for simplicity, we assume that every tt-category~$\cat{K}$ that we discuss below is either rigid or that the phrase `tt-ideal' means `radical tt-ideal'.
\end{Rem}

\begin{Not}
For a class $\cat{E}\subseteq\cat{K}$ of objects, the tt-ideal generated by~$\cat{E}$ is~$\ideal{\cat{E}}=\bigcap_{\cat{J}\supseteq\cat{E}}\cat{J}$, where $\cat{J}$ runs through the tt-ideals containing~$\cat{E}$.
\end{Not}

\begin{Def}
\index{prime tt-ideal}%
\index{spectrum (in tt-geometry)}%
\index{support (of an object)}%
A \emph{prime} $\cat{P}\subset\cat{K}$ is a proper tt-ideal such that $X\otimes Y\in\cat{P}$ forces $X\in\cat{P}$ or $Y\in\cat{P}$. We denote the set of prime tt-ideals by
\[
\SpcK=\SET{\cat{P}\subset\cat{K}}{\cat{P}\textrm{ is prime}}
\]
and call it the \emph{spectrum} of~$\cat{K}$. The \emph{support} of an object $X\in\cat{K}$ is the subset
\[
\supp(X)=\SET{\cat{P}\in\SpcK}{X\notin\cat{P}}\,.
\]
The topology of~$\SpcK$ is defined to have $\{\supp(X)\}_{X\in \cat{K}}$ as basis of closed subsets. Explicitly, for each set of objects $\cat{E}\subseteq\cat{K}$ the subset~$U(\cat{E})=\SET{\cat{P}\in\SpcK}{\cat{E}\cap \cat{P}\neq\varnothing}$ is an open of~$\SpcK$, and all open subsets are of this form, for some~$\cat{E}$.
\end{Def}

\begin{Rem}
\label{rem:universal}%
The above construction is introduced in~\cite{Balmer05a}, where the pair $(\SpcK,\supp)$ is characterized by a universal property\,: It is the \emph{final support data}. See~\cite[Thm.\,3.2]{Balmer05a}. We shall not make this explicit but intuitively it means that the space~$\SpcK$ is the best possible one carrying \emph{closed} supports for objects of~$\cat{K}$ with the following rules for all $X,Y,Z$ in~$\cat{K}$\,:
\begin{enumerate}[(1)]
\item $\supp(0)$ is empty and $\supp(\unit)$ is the whole space;
\item $\supp(X\oplus Y)=\supp(X)\cup \supp(Y)$;
\item $\supp(\Sigma X)=\supp(X)$;
\item $\supp(Z)\subseteq\supp(X)\cup \supp(Y)$ for each exact triangle $X\!\to \!Y\!\to \!Z\to \!\Sigma X$;
\item $\supp(X\otimes Y)=\supp(X)\cap \supp(Y)$.
\end{enumerate}
\end{Rem}

\begin{Exa}
\label{exa:P-bar}%
Dually to the Zariski topology, in the tt-spectrum~$\SpcK$ the closure of a point~$\cat{P}\in\SpcK$ consists of all the primes \emph{contained} in it: $\adhpt{\cat{P}}=\SET{\cat{Q}\in\SpcK}{\cat{Q}\subseteq\cat{P}}$. See~\cite[Prop.\,2.9]{Balmer05a}.
\end{Exa}

\begin{Rem}
\label{rem:spectral}%
The tt-spectrum~$\SpcK$ is always a \emph{spectral space} in the sense of Hochster~\cite{Hochster69}\,: it is quasi-compact, it admits a basis of quasi-compact open subsets, and each  of its non-empty irreducible closed subsets has a unique generic point. See~\cite[\S\,2]{Balmer05a}.
\end{Rem}

\begin{Rem}
\label{rem:functorial}%
The construction~$\cat{K}\mapsto \SpcK$ is a contravariant functor. Every exact $\otimes$-functor $F\colon \cat{K}\to \cat{K}'$ between tt-categories induces a continuous (spectral) map $\varphi=\Spc(F)\colon \Spc(\cat{K}')\to \SpcK$ defined by~$\cat{Q}\mapsto F^{-1}(\cat{Q})$. It satisfies $\varphi^{-1}(\supp(X))=\supp(F(X))$ for all~$X\in\cat{K}$. See~\cite[\S\,3]{Balmer05a}.
\end{Rem}

To express tt-classification via the spectrum, we need some preparation.
\begin{Def}
\label{def:K_V}%
To every subset $V\subseteq \SpcK$ we can associate a tt-ideal
\[
\cat{K}_V=\SET{X\in \cat{K}}{\supp(X)\subseteq V}
\]
of~$\cat{K}$. (In fact, this tt-ideal is always radical. See Remark~\ref{rem:rigid}.)
\end{Def}

\begin{Def}
A subset $V\subseteq\SpcK$ is called a \emph{Thomason subset} if it is the union of the complements of a collection of quasi-compact open subsets: $V=\cup_{\alpha}V_\alpha$ where each $V_{\alpha}$ is closed with quasi-compact complement. In the terminology of Hochster~\cite{Hochster69}, these are the \emph{dual-open} subsets.
\end{Def}

\begin{Exa}
\label{exa:noetherian}%
If the space~$\SpcK$ is topologically noetherian (\ie all open subsets are quasi-compact), then $V$ being Thomason is just being closed under specialization ($x\in V\Rightarrow\adhpt{x}\subseteq V$), \ie being a union of closed subsets.
\end{Exa}

\begin{Thm}[Classification of tt-ideals, {\cite[Thm.\,4.10]{Balmer05a}}]
\label{thm:classification}%
The assignment $V\mapsto \cat{K}_V$ of Definition~\ref{def:K_V} defines an order-preserving bijection between the Thomason subsets $V\subseteq\SpcK$ and the (radical) tt-ideals~$\cat{J}\subseteq\cat{K}$ of~$\cat{K}$, whose inverse is given by~$\cat{J}\mapsto\supp(\cat{J}):=\cup_{X\in\cat{J}}\supp(X)=\SET{\cat{P}}{\cat{J}\not\subseteq\cat{P}}$.
\end{Thm}

Specifically for the tt-classification of objects $X,Y\in\cat{K}$ (see Remark~\ref{rem:rigid})\,:
\begin{Cor}
\label{cor:classification}%
Two objects~$X,Y\in\cat{K}$ generate the same tt-ideals $\ideal{X}=\ideal{Y}$ if and only if they have the same support~$\supp(X)=\supp(Y)$. More precisely, $Y$ belongs to~$\ideal{X}$ if and only if~$\supp(Y)\subseteq\supp(X)$.
\end{Cor}

The following converse to Theorem~\ref{thm:classification} holds. See~\cite{Balmer05a} for details.
\begin{Thm}[Balmer/Buan-Krause-Solberg]
\label{thm:computation}%
Suppose that a spectral space~$\mathcal{S}$ carries a support data $\sigma(X)\subseteq \mathcal{S}$ for~$X\in \cat{K}$ in the sense of~\cite{Balmer05a} and suppose that the assignment $\mathcal{S}\supseteq V\mapsto \SET{X\in\cat{K}}{\sigma(X)\subseteq V}$ induces a bijection between Thomason subsets~$V$ of~$\mathcal{S}$ and (radical) tt-ideals of~$\cat{K}$. Then the canonical map $\mathcal{S}\to \SpcK$ of Remark~\ref{rem:universal} is a homeomorphism.
\end{Thm}

This result was established in~\cite{Balmer05a} under the additional assumption that~$\mathcal{S}$ be noetherian. It was proved in the above maximal generality in~\cite{BuanKrauseSolberg07}. (See Remark~\ref{rem:spectral}.)

\begin{Rem}
Theorems~\ref{thm:classification} and~\ref{thm:computation} allow for a compact reformulation of tt-classifications, including the ones anterior to~\cite{Balmer05a}. Thus most classifications for the tt-categories~$\cat{K}$ discussed in Sections~\ref{sec:topology}-\ref{sec:motives} are phrased in the simple form of a description of~$\SpcK$. The tt-classification is then always the same, in terms of subsets of~$\SpcK$, as in Theorems~\ref{thm:classification} and Corollary~\ref{cor:classification}, \emph{and we shall not repeat these corollaries}.

On the other hand, approaching tt-classification via~$\SpcK$ buys us some flexibility, for partial results about $\SpcK$ can be interesting while a `partial classification' is an odd concept. For instance, one can know $\SpcK$ \emph{as a set} in some examples, with partial information on the topology. Or one can describe $\SpcK=U\cup Z$ with a complete description of the closed subset~$Z$ and its open complement~$U$ without knowing exactly how they attach. And so on.

In recent years, the geometric study of the tt-spectrum \textsl{per se} has led to new computations of~$\SpcK$, from which the tt-classification can be deduced a posteriori. This will be illustrated in the later sections.
\end{Rem}

\begin{Rem}
Some of the above results connect to lattice theory, see \cite{BuanKrauseSolberg07,KockPitsch17}. It is a non-trivial property of a lattice, like that of tt-ideals in~$\cat{K}$, to be \emph{spatial}, \ie in bijection with the open subsets of a topological space. In fact, without the tensor this fails in general (Remark~\ref{rem:no-tensor-no-cry}).

The tt-classification of Theorem~\ref{thm:classification} tacitly assumes that $\cat{K}$ consists of `small enough' objects. Assumption~\ref{ass:small} and Remark~\ref{rem:rigid} belong to this logic too. Another indication of the smallness of~$\cat{K}$ is that we do not mention infinite coproducts in~$\cat{K}$, and we only discuss thick subcategories, not localizing ones (\ie those closed under arbitrary coproducts). When dealing with a `big' tt-category~$\cat{T}$, the natural candidate for a `small'~$\cat{K}$ is the subcategory of rigid objects in~$\cat{T}$, which may or may not coincide with compact ones.

There are also `big' subcategories of `big' tt-categories worth investigating, most famously \emph{smashing} subcategories. It is an open problem whether the lattice of smashing $\otimes$-ideals is spatial or not. We prove in~\cite{BalmerKrauseStevenson17app} that it is a \emph{frame}, thus it is at least `spatial' in the quirky sense of pointless topology.

The connection between thick subcategories of compact objects and smashing subcategories is a topic in its own right, often dubbed the \emph{Telescope Conjecture}. We shall not attempt to discuss it systematically here but will mention it in a few examples. See Krause~\cite{Krause00} for a beautiful abstract answer via ideals of morphisms.
\end{Rem}

%----------------------------------------------------------------------
\section{Topology}
\label{sec:topology}%
%----------------------------------------------------------------------

As already said, tt-classification (or at least `t-classification') was born in topology, more precisely in \emph{chromatic homotopy theory}, see~\cite{BarthelBeaudry19}. The tt-category we consider here is the topological stable homotopy category~$\SH$, \ie the homotopy category of topological spectra, and more specifically its subcategory~$\SHc$ of compact objects. See for instance~\cite{Ravenel92}. In other words, $\SHc$ is the Spanier-Whitehead stable homotopy category of finite pointed CW-complexes.

The first operation one can do on~$\SH$ is to $p$-localize it at a prime~$p$, \ie invert multiplication by every prime different from~$p$. On compacts, this gives us~$\SHcp$. Both $\SHc$ and~$\SHcp$ are essentially small rigid tt-categories.

\begin{Rem}
\label{rem:no-tensor-in-SH}%
Something special happens in~$\SHc$ and therefore in~$\SHcp$ as well: The unit~$\unit=S^0$, \aka the sphere spectrum, generates the category as a thick triangulated subcategory. Consequently, every thick subcategory is automatically a tt-ideal. In such situations, the tensor is not essential in the tt-classification and we are equivalently classifying thick subcategories.
\end{Rem}

\begin{Rem}
\label{rem:Morava}%
\index{Morava K-theories}%
A critical ingredient in chromatic theory is the countable family of so-called \emph{Morava $K$-theories}, which are homology theories~$K_{p,n}$, for~$n\ge1$, defined on~$\SHcp$ and taking values in graded modules over the `graded field'~$\bbF_p[v_n,v_n^{-1}]$, with $v_n$ in degree~$2(p^n-1)$. See~\cite[\S\,1.5]{Ravenel92}.
\end{Rem}

\begin{Thm}[Hopkins-Smith~{\cite{HopkinsSmith98}}]
\label{thm:HopkinsSmith}%
\index{chromatic classification}%
The spectrum of the classical stable homotopy category~$\SHc$ is the following topological space\,:
\[
\xymatrix@C=.8em @R=.4em{
&&\cat{P}_{2,\infty} \ar@{-}[d]
&\cat{P}_{3,\infty} \ar@{-}[d]
&& \kern-2em{\cdots}
&\cat{P}_{p,\infty} \ar@{-}[d]
& {\cdots}
\\
\Spc(\SHc)=
&&{\vdots} \ar@{-}[d]
& {\vdots} \ar@{-}[d]
&&& {\vdots} \ar@{-}[d]
\\
&&\cat{P}_{2,n+1} \ar@{-}[d]
& \cat{P}_{3,n+1} \ar@{-}[d]
&& \kern-2em{\cdots}
& \cat{P}_{p,n+1} \ar@{-}[d]
& {\cdots}
\\
&&\cat{P}_{2,n} \ar@{-}[d]
& \cat{P}_{3,n} \ar@{-}[d]
&& \kern-2em{\cdots}
& \cat{P}_{p,n} \ar@{-}[d]
& {\cdots}
\\
&&{\vdots} \ar@{-}[d]
& {\vdots} \ar@{-}[d]
&&& {\vdots} \ar@{-}[d]
\\
&&\cat{P}_{2,2} \ar@{-}[rrd]
& \cat{P}_{3,2} \ar@{-}[rd]
&& \kern-2em{\cdots}
& \cat{P}_{p,2} \ar@{-}[lld]
& {\cdots}
\\
&&&& \cat{P}_{0,1}
}
\]
in which every line indicates that the higher point belongs to the closure of the lower one (Example~\ref{exa:P-bar}). More precisely:
\begin{enumerate}[\rm(a)]
\item The tt-prime $\cat{P}_{0,1}$ is the kernel of rationalization $\SHc\to \SHc_\bbQ\cong\Db(\bbQ)$, that is, the subcategory of torsion spectra. It is the dense point of~$\Spc(\SHc)$.
\item For each prime number~$p$, the tt-prime $\cat{P}_{p,\infty}$ is the kernel of localization $\SHc\to \SHcp$. These $\cat{P}_{p,\infty}$ are exactly the closed points of~$\Spc(\SHc)$.
\item For each prime number~$p$ and each integer~$2\le n<\infty$, the tt-prime~$\cat{P}_{p,n}$ is the kernel of the composite $\SHc\to \SHcp \to \bbF_p[v_{n-1}^{\pm1}]\grmod$ of localization at~$p$ and $(n-1)^\textrm{st}$ Morava $K$-theory~$K_{p,n-1}$ (Remark~\ref{rem:Morava}).
\item The support of an object~$X$ is either empty when $X=0$, or the whole of~$\Spc(\SHc)$ when $X$ is non-torsion, or a finite union of `columns'
    \[
    \adhpt{\cat{P}_{p,m_p}}=\SET{\cat{P}_{p,n}}{m_p\le n\le\infty}
    \]
    for integers $2\le m_p<\infty$ varying with the primes~$p$.
\item A closed subset is either empty, or the whole $\Spc(\SHc)$, or a finite union of closed points~$\{\cat{P}_{p,\infty}\}$ and of columns $\adhpt{\cat{P}_{p,m_p}}$ with $m_p\ge 2$ as in~(d).
\item A Thomason subset of~$\Spc(\SHc)$ is either empty, or the whole $\Spc(\SHc)$, or an arbitrary union of columns $\adhpt{\cat{P}_{p,m_p}}$ with $m_p\ge2$ as in~(d).
\end{enumerate}
\end{Thm}

The above Theorem~\ref{thm:HopkinsSmith} is not the way the chromatic filtration is expressed in the original literature; see the translation in~\cite[\S\,9]{Balmer10b}.

\begin{Exa}
\label{exa:p-torsion}%
An object~$X\in\SHc$ has support contained in the $p$-th column, $\supp(X)\subseteq\adhpt{\cat{P}_{p,2}}$, if and only if it is `$p$-primary torsion', \ie it satisfies $p^\ell\cdot X=0$ for some $\ell\ge 1$.
\end{Exa}

\begin{Exa}
\label{exa:torsion}%
The support of the tt-ideal $\cat{J}=\cat{P}_{0,1}$ of torsion spectra is exactly the Thomason subset $\Spc(\SHc)\setminus\{\cat{P}_{0,1}\}$ and is therefore the disjoint union of all columns $\sqcup_{p}\adhpt{\cat{P}_{p,2}}$. This reflects the fact that a torsion object in~$\SHc$ is the direct sum of $p$-primary torsion objects as in Example~\ref{exa:p-torsion}.
\end{Exa}

\begin{Rem}
The fact that the closed point~$\{\cat{P}_{p,\infty}\}$ cannot be the support of an object reflects the fact that an object in~$\SHcp$ which is killed by all Morava $K$-theories~$K_{p,n}$ for $n\ge1$ must be zero. It also shows that $\SpcK$ is not noetherian, already in this initial case of~$\cat{K}=\SHc$ (see Example~\ref{exa:noetherian}).
\end{Rem}

\begin{Rem}
In this setting, the Telescope Conjecture is open (again). See~\cite{Krause00} and further references therein.
\end{Rem}

%----------------------------------------------------------------------
\section{Commutative algebra and algebraic geometry}
\label{sec:AG}%
%----------------------------------------------------------------------

\index{perfect complexes}%

As already indicated, Hopkins~\cite{Hopkins87} initiated the transposition of the chromatic classification from topology to algebra. The correct statement for noetherian rings was proved by Neeman~\cite{Neeman92a} and the perfect version for general schemes, not necessarily noetherian, is due to Thomason in his last published paper~\cite{Thomason97}. In terms of tt-spectra it becomes the following very beautiful result.
\begin{Thm}[Thomason~{\cite{Thomason97}}]
\label{thm:Thomason}%
Let $\mathcal{X}$ be a scheme which is quasi-compact and quasi-separated. Then the spectrum of the derived category~$\Dperf(\mathcal{X})$ of perfect complexes (with $\otimes={}^{\rmL}\otimes_{\mathcal{O}_{\mathcal{X}}}$) is isomorphic to the underlying space~$|\mathcal{X}|$ itself, via the homeomorphism
\[
\xymatrix@R=.3em{
|\mathcal{X}| \ar[r]^-{\cong}
& \Spc(\Dperf(\mathcal{X}))
\\
x \ar@{|->}[r]
& \cat{P}(x)
}
\]
where, for each point~$x$ of~$\mathcal{X}$, the tt-prime $\cat{P}(x)=\SET{Y\in\Dperf(\mathcal{X})}{Y_{x}\cong 0}$ is the kernel of localization~$\Dperf(\mathcal{X})\to\Dperf(\mathcal{O}_{\mathcal{X},x})$ at~$x$.
\end{Thm}

\begin{Rem}
Equivalently, $\cat{P}(x)$ can be described as the kernel of the residue functor $\Dperf(\mathcal{X})\to \Db(\kappa(x))$ to the residue field~$\kappa(x)$ of~$\mathcal{X}$ at~$x$.
\end{Rem}

\begin{Rem}
Recall that a scheme~$\mathcal{X}$ is quasi-compact and quasi-separated if the underlying space~$|\mathcal{X}|$ admits a basis of quasi-compact open subsets (including~$|\mathcal{X}|$ itself). This purely topological condition is equivalent to $|\mathcal{X}|$ being spectral. Hence this condition is the maximal generality in which the above result can hold in view of Remark~\ref{rem:spectral}. Noetherian schemes and affine schemes are quasi-compact and quasi-separated.
\end{Rem}

The affine case vindicates our use of the word `spectrum'\,(\footnote{\,If this creates confusion with Example~\ref{exa:P-bar}, note that the map of Corollary~\ref{cor:Thomason}, $\gp\mapsto \cat{P}(\gp)=\Ker(\Dperf(A)\to \Dperf(A_\gp))$, \emph{reverses} inclusions: $\gp\subseteq\gq\Rightarrow \cat{P}(\gp)\supseteq\cat{P}(\gq)$.}):
\begin{Cor}
\label{cor:Thomason}%
Let $A$ be a commutative ring. Then the tt-spectrum of the homotopy category~$\Kb(A\pproj)\cong\Dperf(A)$ of bounded complexes of finitely generated projective~$A$-modules is homeomorphic to the Zariski spectrum of~$A$
\[
\Spec(A)\isoto\Spc(\Kb(A\pproj))\,.
\]
\end{Cor}

\begin{Rem}
An error in \cite{Hopkins87}, corrected in~\cite{Neeman92a}, was not to assume~$A$ noetherian. However we see that Thomason's Corollary~\ref{cor:Thomason} does not assume $A$ noetherian. The point is that the tt-classification (Theorem~\ref{thm:classification}) which is equivalent to Corollary~\ref{cor:Thomason} involves actual Thomason subsets not mere specialization-closed subsets, whereas~\cite{Hopkins87} and~\cite{Neeman92a} are phrased in terms of specialization-closed subsets. The assumption that $A$ is noetherian is only useful to replace `Thomason' by `specialization-closed' (Example~\ref{exa:noetherian}).
\end{Rem}

\begin{Rem}
\label{rem:no-tensor-no-cry}%
As we saw in the topological example of Section~\ref{sec:topology}, when the unit~$\unit$ generates the tt-category~$\cat{K}$ as a thick subcategory we do not really need the tensor. This is also the case for~$\cat{K}=\Kb(A\pproj)$ for instance.

But in general the tensor is essential for classification by means of subsets of~$\SpcK$. Indeed, the lattice of thick subcategories of a triangulated category $\cat{K}$ cannot be classified in terms of the lattice of subsets of pretty much anything because it may not satisfy \emph{distributivity}: $\cat{J}_1\wedge(\cat{J}_2\vee\cat{J}_3)=(\cat{J}_1\wedge\cat{J}_2)\vee(\cat{J}_1\wedge\cat{J}_3)$. Already for $\cat{K}=\Dperf(\mathcal{X})$ over the projective line~$\mathcal{X}=\bbP^1_{k}$ distri\-bu\-tivity fails with $\cat{J}_i$ the thick subcategory generated by~$\mathcal{O}(i)$. See~\cite[Rem.\,5.10]{BalmerKrauseStevenson17app}.
\end{Rem}

\begin{Rem}
An application of Theorem~\ref{thm:Thomason} is the reconstruction of every quasi-compact and quasi-separated scheme~$\mathcal{X}$ from the data of the \emph{tensor}-triangulated category~$\Dperf(\mathcal{X})$. Indeed, one can equip the tt-spectrum~$\SpcK$ with a sheaf of commutative rings, which in the case of~$\cat{K}=\Dperf(\mathcal{X})$ recovers the structure sheaf~$\mathcal{O}_{\mathcal{X}}$. See details in~\cite[\S\,6]{Balmer05a}. By contrast, Mukai~\cite{Mukai81} proved earlier that such a reconstruction is impossible from the triangulated structure alone.
\end{Rem}

\begin{Rem}
The above~$\cat{K}=\Dperf(R)$ are the compact and rigid objects in the big derived category~$\cat{T}=\Der(R)$. The Telescope Conjecture holds when $R$ is noetherian by Neeman~\cite{Neeman92a} but fails in general by Keller~\cite{Keller94b}.
\end{Rem}

\begin{Rem}
Other tt-categories can be associated to schemes, or commutative rings, for instance right-bounded derived categories. For first results in this direction, see work of Matsui and Takahashi~\cite{MatsuiTakahashi17,Matsui18}.
\end{Rem}

One can generalize Theorem~\ref{thm:Thomason} almost verbatim to reasonable stacks:

\begin{Thm}[Hall~{\cite[Thm.\,1.2]{Hall16}}]
Let $\mathcal{X}$ be a quasi-compact algebraic stack with quasi-finite separated diagonal, whose stabilizer groups at geometric points are finite
linearly reductive group schemes ($\mathcal{X}$ is `tame'). Then $\Spc(\Dperf(\mathcal{X}))\cong|\mathcal{X}|$.
\end{Thm}

We refer to~\cite{Hall16} for terminology. Note earlier work of Krishna~\cite{Krishna09} in characteristic~0, and of Dubey-Mallick~\cite{DubeyMallick12} for finite groups acting on smooth schemes in characteristic prime to the order of the groups.

\bigbreak

One can also consider the graded version of Corollary~\ref{cor:Thomason}:
\begin{Thm}[Dell'Ambrogio-Stevenson~{\cite[Thm.\,4.7]{DellAmbrogioStevenson14}}]
Let $A$ be a graded-commutative ring (graded over any abelian group), then there is a canonical isomorphism $\Spc(\Dperf(A))\cong\Spech(A)$, between the tt-spectrum of~$\Dperf(A)$ and the spectrum of homogeneous prime ideals of~$A$.
\end{Thm}

Let us mention a variation relating to singularities.
\index{singularity categories}%
\begin{Thm}[Stevenson~{\cite[Thm.\,7.7]{Stevenson14}}]
\label{thm:singularity}%
Let $\mathcal{X}$ be a noetherian separated scheme with only hypersurface singularities. Then there is an order-preserving bijection between the specialization-closed subsets of the singular locus of~$\mathcal{X}$ and the thick $\Dperf(\mathcal{X})$-submodules of the singularity category~$\Db(\coh \mathcal{X})/\Dperf(\mathcal{X})$.
\end{Thm}

Here the singularity category is not itself a tt-category but a triangulated category with an action by the tt-category~$\Dperf(\mathcal{X})$. As such, this result is an application of Stevenson's \emph{relative tt-geometry}~\cite{Stevenson13}. Another application of Stevenson's theory is the tt-classification for derived categories of matrix factorizations in Hirano~\cite{Hirano19}, which extends earlier result of Takahashi~\cite{Takahashi10}.

%----------------------------------------------------------------------
\section{Modular representation theory and related topics}
\label{sec:modular}%
%----------------------------------------------------------------------

%
\begin{Point}
\index{stable module category}%
Let $G$ be a finite group and let $k$ be a field. Maschke's Theorem says that the order of~$G$ is invertible in~$k$ if and only if $kG$ is semisimple. In that case, all $kG$-modules are projective. \emph{Modular} representation theory refers to the non-semisimple situation. Then the \emph{stable module category} is the additive quotient~\cite{Happel88}
\[
kG\stmod=\frac{kG\mmod}{kG\pproj}
\]
which precisely measures how far $kG$ is from being semisimple. It is a tt-category whose objects are all finitely generated $kG$-modules and whose groups of morphisms $\Hom_{kG\stmod}(X,Y)$ are given by the quotient of the abelian group of $kG$-linear maps $\Hom_{kG}(X,Y)$ modulo the subgroup of those maps factoring via a projective module. Tensor is over~$k$ with diagonal $G$-action: $g\cdot (x\otimes y)=(gx)\otimes (gy)$ in~$X\otimes_k Y$. The $\otimes$-unit is $\unit=k$ with trivial $G$-action.
\end{Point}

\begin{Point}
We can also consider the derived category $\Db(kG\mmod)$, with the `same' tensor. Every non-zero tt-ideal $\cat{J}\subseteq\Db(kG\mmod)$ contains $\Dperf(kG)$ because $kG\otimes-\cong \Ind_{1}^G\Res^G_{1}$ and $\Db(k\mmod)=\Dperf(k)$ is semisimple. Hence the tt-classification of~$\Db(kG\mmod)$ and of its Verdier quotient by~$\Dperf(kG)$ are very close. (The former has just one more tt-ideal: zero.) By Rickard~\cite{Rickard89}, that quotient is equivalent to the stable module category:
\[
\frac{\Db(kG\mmod)}{\Dperf(kG)}\cong kG\stmod.
\]
\end{Point}

\begin{Thm}[Benson-Carlson-Rickard~{\cite{BensonCarlsonRickard97}}]
\label{thm:finite-group}%
\index{support variety (of a finite group)}%
There is a homeomorphism between the spectrum of the stable module category and the so-called projective support variety
\[
\Spc(kG\stmod)\cong \Proj(\rmH^\sbull(G,k))
\]
which can be extended (by adding one closed point) to a homeomorphism
\[
\Spc(\Db(kG\mmod))\cong \Spech(\rmH^\sbull(G,k)).
\]
Explicitly, to every homogeneous prime~$\gp^\sbull\subset\rmH^\sbull(G,k)$ corresponds the tt-prime
\[
\cat{P}(\gp^\sbull)=\SET{X}{\textrm{ there is a homogeneous }\zeta\notin\gp^\sbull\textrm{ such that }\zeta\cdot X=0}\,.
\]
\end{Thm}

Again, the above does not appear verbatim in the source. See details in~\cite{Balmer05a} or~\cite[Prop.\,8.5]{Balmer10b}. A more recent proof can be found in~\cite{CarlsonIyengar15}.

\begin{Rem}
The reader interested in the related derived category of cochains on the classifying space~$BG$ is referred to~\cite{BensonIyengarKrause11b} for finite groups and to the comprehensive recent work~\cite{BarthelCastellanaHeardValenzuela19} for $p$-local compact groups; see comments and references therein about earlier work by Benson-Greenlees.
\end{Rem}

\begin{Point}
For finite group~\emph{schemes}~$G$, the following generalization of Theorem~\ref{thm:finite-group} would follow from claims made in~\cite{FriedlanderPevtsova07} but a flaw was found in~\cite[Thm.\,5.3]{FriedlanderPevtsova07}, which was eventually fixed in the recent~\cite{BensonIyengarKrausePevtsova18}; see in particular \cite[Rem.\,5.4 and Thm.\,10.3]{BensonIyengarKrausePevtsova18}.
\end{Point}

\begin{Thm}[Benson-Friedlander-Iyengar-Krause-Pevtsova]
\label{thm:finite-group-scheme}%
For a finite group scheme~$G$ over~$k$, there is a homeomorphism $\Spc(kG\stmod)\cong\Proj(\rmH^\sbull(G,k))$.
\end{Thm}

\begin{Rem}
Stable module categories of finite group schemes over a field are very `noetherian' and several other results are known about the `big' stable module category as well, like the Telescope Conjecture. See details in~\cite{BensonIyengarKrausePevtsova18}. The technique of \emph{stratification} has led to the tt-classification (of small and large subcategories) in several `noetherian enough' derived settings. See the survey in~\cite{BensonIyengarKrause12} and further references in~\cite{BensonIyengarKrausePevtsova18}.
\end{Rem}

\begin{Point}
Extending beyond field coefficients to other rings~$R$, we can consider the \emph{relative stable module category}~$RG\strel$, obtained from the Frobenius exact structure on the exact category of finitely generated~$RG$-modules with $R$-split exact sequences. Already in small Krull dimension, interesting phenomena can be observed, as in the following result.
\end{Point}

\begin{Thm}[Baland-Chirvasitu-Stevenson~{\cite[Thm.\,1.1]{BalandChirvasituStevenson19}}]
Let $S$ be a discrete valuation ring having residue field~$k$ and uniformizing parameter~$t$ and let $R_n=S/t^n$. Let $G$ be a finite group. Then the tt-spectrum of the relative stable module category
\[
\Spc(R_n G\strel)\cong \sqcup_{i=1}^n \Spc(kG\stmod),
\]
is a coproduct of~$n$ copies of the projective support variety of Theorem~\ref{thm:finite-group}.
\end{Thm}

On the topic of singularity categories, let us mention~\cite{Xu14} and its recent generalization (recall that a category is EI if any endomorphism is invertible):
\begin{Thm}[Wang~{\cite[Thm.\,5.2]{Wang19}}]
\index{singularity categories}%
Let $\cat{C}$ be a finite EI category, projective over a field~$k$ and $\Der_{\textrm{sing}}(k\cat{C})=\Db(k\cat{C}\mmod)/\Db(k\cat{C}\pproj)$ its singularity category. Then there is a homeomorphism
\[
\Spc(\Der_{\textrm{sing}}(k\cat{C}))\cong \sqcup_{x\in\cat{C}} \Spc(kG_x\stmod),
\]
where $G_x=\Aut_{\cat{C}}(x)$.
\end{Thm}

\begin{Point}
Antieau-Stevenson~\cite{AntieauStevenson16} consider further derived categories of representations of small categories over commutative noetherian rings. They obtain several interesting classifications, including for localizing subcategories, and in particular for simply laced Dynkin quivers. See also the earlier~\cite{LiuSierra13}.
\end{Point}

\begin{Point}
Let us now turn our attention to stable module categories related to Lie algebras. Boe-Kujawa-Nakano~\cite{BoeKujawaNakano17} prove several results about classical Lie superalgebras. In particular for the general linear Lie superalgebra $\mathfrak{g}=\mathfrak{gl}(m|n)=\mathfrak{g}_{\bar0}\oplus \mathfrak{g}_{\bar1}$ and $\cat{K}=\underline{\cat{F}}$ the stable category of the category~$\cat{F}$ of finite dimensional~$\mathfrak{g}$-modules which admit a compatible action by~$G_{\bar0}$ and are completely reducible as $G_{\bar0}$-modules (where~$\textrm{Lie}\,G_{\bar0}=\mathfrak{g}_{\bar0}$). They prove in~\cite[Thm.\,5.2.2]{BoeKujawaNakano17} that the spectrum~$\Spc(\underline{\cat{F}})$ is homeomorphic to the $N$-homogenous spectrum $N-\Proj(S^\sbull(\mathfrak{f}_{\bar1}))$ where $\mathfrak{f}$ is the detecting subalgebra of~$\mathfrak{g}$ and $N=\textrm{Norm}_{G_{\bar0}}(\mathfrak{f}_{\bar1})$.
\end{Point}

The same authors more recently considered quantum groups:
\begin{Thm}[Boe-Kujawa-Nakano~{\cite[Thm.\,7.6.1]{BoeKujawaNakano17pp}}]
Let $G$ be a complex simple algebraic group over~$\bbC$ with~$\mathfrak{g}=\textrm{\rm Lie}\,G$. Assume that~$\zeta$ is a primitive~$\ell$th root of unity where~$\ell$ is greater than the Coxeter number for~$\mathfrak{g}$. Then the tt-spectrum of the stable module category for the quantum group~$U_\zeta(\mathfrak{g})$ is
\[
\Spc(U_\zeta(\mathfrak{g})\stmod)\cong G-\Proj(\bbC[\mathcal{N}])
\]
where $\mathcal{N}$ is the nullcone, \ie the set of nilpotent elements of~$\mathfrak{g}$.
\end{Thm}

\begin{Exa}
Another example where $\Spc(A\stmod)$ is isomorphic to the variety~$\Proj(H^\sbull(A,k))$ is the algebra $A=k[X_1,\ldots,X_n]/(X_1^\ell,\ldots,X_n^\ell)\rtimes(\bbZ/\ell\bbZ)^{\times n}$ which appears in Pevtsova-Witherspoon~\cite[Thm.\,1.2]{PevtsovaWitherspoon15}.
\end{Exa}

%----------------------------------------------------------------------
\section{Equivariant stable homotopy and $KK$-theory}
\label{sec:equivariant}%
%----------------------------------------------------------------------

\index{equivariant stable homotopy category}%

\begin{Point}
Let $G$ be a compact Lie group, \eg\ a finite group, and let~$\SH(G)$ be the equivariant stable homotopy category of genuine $G$-spectra. The tensor-triangulated category of compact (rigid) objects in~$\SH(G)$ is denoted~$\SH(G)^c$. In general, the spectrum of~$\SH(G)^c$ is not quite known but significant progress occurred in recent years. It relies in an essential way on the non-equivariant case~$G=1$ of Section~\ref{sec:topology}.
\end{Point}

\begin{Point}
\label{pt:P(H,p,n)}%
For every chromatic tt-prime $\cat{P}_{p,n}\in\Spc(\SHc)$ in the stable homotopy category~$\SHc$ (Theorem~\ref{thm:HopkinsSmith}) and every closed subgroup~$H\le G$, let
\[
\cat{P}(H,p,n)=(\Phi^H)^{-1}(\cat{P}_{p,n})
\]
be its preimage under geometric $H$-fixed points $\Phi^{H}\colon \SH(G)^c\to \SHc$, which is a tt-functor. This `equivariant' tt-prime~$\cat{P}(H,p,n)$ is the image of the chromatic $\cat{P}_{p,n}$ under $\Spc(\Phi^H)\colon \Spc(\SHc)\to \Spc(\SH(G)^c)$ as in Remark~\ref{rem:functorial}.

It is convenient to use the convention that $\cat{P}_{p,1}$ means~$\cat{P}_{0,1}$ for all~$p$. And similarly, to read~$\cat{P}(H,p,1)$ as~$\cat{P}(H,0,1)$.
\end{Point}

Let us first discuss the case where $G$ is a finite group. Varying the subgroup~$H\le G$, the maps $\Spc(\Phi^H)$ cover~$\Spc(\SH(G)^c)$ -- a fact that is also true for general compact Lie groups, see Theorem~\ref{thm:SH(G)-Lie}.

\begin{Thm}[Balmer-Sanders~{\cite{BalmerSanders17}}]
\label{thm:SH(G)}%
Let $G$ be a finite group. Then every tt-prime in~$\SH(G)^c$ is of the form~$\cat{P}(H,p,n)$ for a unique subgroup~$H\le G$ up to conjugation and a unique chromatic tt-prime $\cat{P}_{p,n}\in\Spc(\SHc)$. Understanding inclusions between tt-primes completely describes the topology on~$\Spc(\SH(G)^c)$.

If $K\lhd H$ is a normal subgroup of index~$p>0$, then $\cat{P}(K,p,n+1)\subset \cat{P}(H,p,n)$ for every~$n\ge 1$. There is no inclusion $\cat{P}(K,q,n)\subseteq\cat{P}(H,p,m)$ unless the corresponding chromatic tt-primes are included $\cat{P}_{q,n}\subseteq\cat{P}_{p,m}$ (which forces $n\ge m$, and $p=q$ if $m>1$) \emph{and} $K$ is conjugate to a $q$-subnormal subgroup of~$H$ (see~\ref{pt:subnormal}).
\end{Thm}

\begin{Point}
For finite groups of square-free order, like $G=C_p$ for instance, the above result completely describes $\Spc(\SH(G)^c)$, with its topology, and thus gives the tt-classification. This result was a first major example where~$\SpcK$ was determined first and the tt-classification deduced as a corollary.
\end{Point}

\begin{Point}
\label{pt:subnormal}%
For other groups, the question is to decide when $\cat{P}(K,p,n)\subset\cat{P}(H,p,m)$, in terms of~$n-m$, for $K\le H$ a $p$-subnormal subgroup of~$H$ (\ie one such that there exists a tower of normal subgroups of index~$p$ from~$K$ to~$H$). Theorem~\ref{thm:SH(G)} implies that this inclusion holds when $n-m\ge \log_p([H:K])$.
\end{Point}

The case of abelian groups (and a little more) was recently tackled in~\cite{BHNNNS19}, showing that the above~$\log_p([H:K])$ is not the sharpest bound.
\begin{Thm}[Barthel-Hausmann-Naumann-Nikolaus-Noel-Stapleton]
\label{thm:BHNNNS}%
Let $G$ be a finite abelian group, let $K\le H\le G$ be subgroups, let $p$ be a prime and let~$1\le n<\infty$ be an integer. Then the minimal~$i$ such that $\cat{P}(K,p,n)\subseteq\cat{P}(H,p,n-i)$ is
$i= \rk_p(H/K)$ the $p$-rank of the quotient.
\end{Thm}

See~\cite[Cor.\,1.3]{BHNNNS19}. The precise topology of $\Spc(\SH(G)^c)$ for general finite groups remains an open problem.

\begin{Point}
Let us now consider the case of an arbitrary compact Lie group~$G$ but after rationalization~$\SH(G)^c_\bbQ$. Of course, tensoring with~$\bbQ$ hides the `chromatic direction' but this is an essential step in understanding the `equivariant direction', by which we mean the role played by the subgroups of~$G$.

A closed subgroup $K\le H$ is \emph{cotoral} if $K$ is normal and $H/K$ is a torus. Every closed subgroup~$H\le G$ defines a tt-prime~$\cat{P}_H$ in the spectrum~$\Spc(\SH(G)^c_\bbQ)$, namely $\cat{P}_H:=\Ker(\Phi^H)$ the kernel of geometric $H$-fixed points, \ie the image of the unique prime $(0)$ under the map $\Spc(\Db(\bbQ))\to \Spc(\SH(G)^c)$ associated to the tt-functor $\Phi^H\colon \SH(G)^c_\bbQ\to \SH^c_\bbQ\cong\Db(\bbQ)$.
\end{Point}

\begin{Thm}[Greenlees~{\cite[Thm.\,1.3]{Greenlees19}}]
Let $G$ be a compact Lie group. Every tt-prime of the rational equivariant stable homotopy category~$\SH(G)^c_\bbQ$ is equal to~$\cat{P}_H=\Ker(\Phi^H)$ for a closed subgroup~$H\le G$, unique up to conjugation. Furthermore, specialization of tt-primes corresponds to cotoral inclusions: We have $\cat{P}_K\subseteq\cat{P}_H$ if and only if $K$ is conjugate to a cotoral subgroup of~$H$. The topology on~$\Spc(\SH(G)^c_\bbQ)$ corresponds to the $f$-topology of~\cite{Greenlees98}.
\end{Thm}

Recently there has been further progress for arbitrary compact Lie group:
\begin{Thm}[Barthel-Greenlees-Hausmann~{\cite{BarthelGreenleesHausmann18pp}}]
\label{thm:SH(G)-Lie}%
Let $G$ be a compact Lie group. Then every tt-prime of $\SH(G)^c$ is of the form $\cat{P}(H,p,n)$ as in~\ref{pt:P(H,p,n)}. Moreover, the topology is completely understood in terms of inclusions of tt-primes.
\end{Thm}

Barthel-Greenlees-Hausmann more precisely track the inclusion of primes, in terms of functions on the compact and totally-disconnected Hausdorff orbit space $\textrm{Sub}(G)/G$ of $G$ acting by conjugation on its closed subgroups. Furthermore, they give a complete description of the topology in the case of an \emph{abelian} compact Lie group, extending Theorem~\ref{thm:BHNNNS}; see~\cite[Thm.\,1.4]{BarthelGreenleesHausmann18pp}.

\begin{center}
*\ *\ *
\end{center}

\begin{Point}
\label{pt:KK}%
\index{Kasparov's KK-theory}%
The closest to analysis that tt-geometry has gone so far is in the theory of $C^*$-algebras, via Kasparov's $KK$-theory. Although this is not strictly speaking equivariant homotopy theory, we include it in this section as $KK$-theory belongs to the broad topic of noncommutative topology.

One begins with the `cellular' subcategory, \aka the `bootstrap' category.
\end{Point}

\begin{Thm}[Dell'Ambrogio~{\cite[\S\,6]{DellAmbrogio10}}]
Let $G$ be a finite group and $\cat{K}^G$ be the thick subcategory of the $G$-equivariant Kasparov category~$KK^G$ generated by the unit. Then the comparison map $\rho\colon\Spc(\cat{K}^G)\to \Spec(R(G))$ to the spectrum of the complex representation ring~$R(G)$, as in Definition~\ref{def:comparison-map}, is surjective and admits a continuous section. In the non-equivariant case, the above map induces a homeomorphism between the tt-spectrum of the so-called `bootstrap category' $\cat{K}^1=\textrm{\rm Boot}$ and~$\Spec(\bbZ)$.
\end{Thm}

\begin{Point}
Dell'Ambrogio conjectures that $\rho\colon\Spc(\cat{K}^G)\to \Spec(R(G))$ is a homeomorphism for all finite groups. The tt-classification for larger $KK$-categories, or for infinite groups, is another interesting open problem.
\end{Point}

%----------------------------------------------------------------------
\section{Motives and $\Aone$-homotopy}
\label{sec:motives}%
%----------------------------------------------------------------------

%
\begin{Point}
We consider here two classes of `motivic' tensor-triangulated categories: First, we have the derived category of motives~$\DM(F;R)$, over a base field~$F$ and with coefficients in a ring~$R$, as first introduced by Voevodsky~\cite{Voevodsky00}; see also~\cite{Ayoub14,CisinskiDeglise09pp}. Secondly, we have $\SH^{\Aone}\!(F)$ the stable $\Aone$-homotopy category over the base~$F$ introduced by Morel and Voevodsky~\cite{MorelVoevodsky99,Voevodsky98,Morel06}. (Other base schemes can be considered.)
\end{Point}

\begin{Point}
\index{Tate motives}%
Let us begin with~$\DM(F;R)$. Like in $KK$-theory (see~\ref{pt:KK}), one can first consider the `cellular' tt-subcategory of \emph{(mixed) Tate motives}~$\DTM(F;R)$ generated as a localizing  subcategory by the invertible Tate objects~$R(i)$ for $i\in\bbZ$. Its subcategory of compact objects is the rigid tt-category~$\DTM(F;R)^c$ we shall discuss now.
\end{Point}

\begin{Point}
\label{pt:Peter}%
Peter~\cite{Peter13} established the tt-bridgehead into motivic territory, when he proved that the spectrum $\Spc(\DTM(F;\bbQ)^c)=\ast$ reduces to a point, for $F$ a field satisfying the Beilinson-Soul\'e vanishing conjecture and a less standard restriction on rational motivic cohomology, namely $\rmH^i_{\mot}(F;\bbQ(j))=0$ for $j\ge i\ge2$. For instance, this applies to~$F=\bar\bbQ$. In fact, $\Spc(\DTM(F;\bbQ)^c)=\ast$ would follow from $\Spc(\DM(F;\bbQ)^c)=\ast$, a conjecture which is supported by:
\end{Point}

\begin{Thm}[Kelly~{\cite[Thm.\,36]{Kelly16pp}}]
\label{thm:Kelly}%
Let $F$ be a finite field such that every connected smooth projective variety~$X$ over~$F$ satisfies the Beilinson-Parshin conjecture and agreement of rational and numerical equivalence. Then $\Spc(\DM(F;\bbQ)^c)=\ast$ is a point.
\end{Thm}

The above are rational results. Our understanding of the integral picture recently evolved thanks to the following breakthrough.
\begin{Thm}[Gallauer~{\cite[Thm.\,8.6]{Gallauer20}}]
\label{thm:Gallauer}%
Let $F$ be an algebraically closed field of characteristic zero\,(\,\footnote{\,In positive characteristic~$\ell$, replace the coefficients~$\bbZ$ by~$\bbZ[1/\ell]$ and only allow~$p\neq\ell$.}) whose rational motivic cohomology~$\rmH^i_{\mot}(F;\bbQ(j))$ vanishes for $i\le0<j$ (Beilinson-Soul\'e) and for $j\ge i\ge2$. Then the spectrum of~$\DTM(F;\bbZ)^c$ is the following\,:
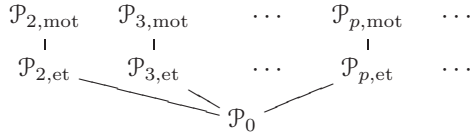
\[
\xymatrix@C=.8em @R=.4em{
&&\cat{P}_{2,\mot} \ar@{-}[d]
& \cat{P}_{3,\mot} \ar@{-}[d]
&& \kern-2em{\cdots}
& \cat{P}_{p,\mot} \ar@{-}[d]
& {\cdots}
\\
&& \cat{P}_{2,\et} \ar@{-}[rrd]
& \cat{P}_{3,\et} \ar@{-}[rd]
&& \kern-2em{\cdots}
& \cat{P}_{p,\et} \ar@{-}[lld]
& {\cdots}
\\
&&&& \cat{P}_{0}
}
\]
where $\cat{P}_0=\Ker(\DTM(F;\bbZ)^c\to \DTM(F;\bbQ)^c)$ consists of the torsion objects and, for every prime number~$p$, the tt-primes $\cat{P}_{p,\mot}$ and $\cat{P}_{p,\et}$ are the kernels of motivic and \'etale cohomology with $\bbZ/p$ coefficients, respectively.
\end{Thm}

\begin{Exa}
Peter's conditions~\ref{pt:Peter} on rational motivic cohomology are satisfied for~$F=\bar{\bbQ}$ for instance. This also provides an example for Theorem~\ref{thm:Gallauer}. Indeed, Gallauer's result uses Peter's theorem rationally.
\end{Exa}

\begin{Rem}
Gallauer~\cite{Gallauer20} proves more general results about the derived category~$\DTM(F;\bbZ/p)^c$ of Tate motives with $\bbZ/p$ coefficients without any condition about rational motivic cohomology. The latter theorem follows from the study of the derived categories of filtered modules in~\cite{Gallauer18}.
\end{Rem}

\begin{Point}
In the case of Theorem~\ref{thm:Kelly}, because $-1$ is a sum of squares in~$F$, the rational derived category of motives coincides with the rational stable $\Aone$-homotopy category $\DM(F;\bbQ)\cong \SH^\Aone(F;\bbQ)$. Let us now mention some integral information about $\SH^\Aone(F)$. The first partial results about its spectrum were obtained in~\cite[\S\,10]{Balmer10b}. The most advanced information is currently:
\end{Point}

\begin{Thm}[Heller-Ormsby~{\cite[Thm.\,1.1]{HellerOrmsby18}}]
Let $F$ be a field of characteristic different from~$2$. Then the comparison map (Definition~\ref{def:comparison-map}) to the homogeneous spectrum of Milnor-Witt $K$-theory is surjective:
\[
\Spc(\SH^\Aone(F)^c)\twoheadrightarrow\Spech(K^{MW}_*(F)).
\]
\end{Thm}

The exact computation of the tt-spectrum of~$\SH(F)^c$ is a major open challenge, which involves understanding the fibers of the above map.

\begin{Rem}
Partial results have also been obtained by Dell'Ambrogio and Tabu\-ada~\cite{DellAmbrogioTabuada12} for non-commutative (dg-)motives.
\end{Rem}

%----------------------------------------------------------------------
\section{Pointers to tt-geometry}
\label{sec:tt-geometry}%
%----------------------------------------------------------------------

\begin{Point}
A snapshot of tensor-triangular geometry as of the year 2010 can be found in~\cite{BalmerICM}. For a more recent survey, see~\cite{Stevenson18}. Beyond those references, let us simply highlight some aspects close to the author's own research.
\end{Point}

A very useful basic tool introduced in~\cite{Balmer10b} is the following \emph{comparison map} between tt-spectra and Zariski spectra of suitable graded rings:

\begin{Def}
\label{def:comparison-map}%
\index{comparison map in tt-geometry}%
Let $u\in\cat{K}$ be a $\otimes$-invertible object and~$R^\sbull_{\cat{K},u}=\oplus_{n\in\bbZ}\Hom_{\cat{K}}(\unit,u\potimes{n})$ the associated graded-commutative graded ring. Then
\[
\cat{P}\mapsto \SET{f\in R^\sbull_{\cat{K},u}}{\cone(f)\notin\cat{P}}
\]
defines a continuous map $\rho^\sbull\colon\SpcK\to \Spech(R^\sbull_{\cat{K},u})$. Without grading, one can similarly define $\rho\colon \SpcK\to \Spec(\Endcat{K}(\unit))$.
\end{Def}

\begin{Point}
Dell'Ambrogio-Stanley~\cite{DellAmbrogioStanley16} give a class of cellular tt-categories for which~$\rho^\sbull$ is a homeomorphism, namely when the ring $R_{\cat{K},\unit}^\sbull$ is concentrated in even degrees and is `regular' in a weak sense, which includes noetherian.
\end{Point}

\begin{Point}
The comparison map was generalized in two directions. First by Dell'\break Ambrogio-Stevenson~\cite{DellAmbrogioStevenson14}, by allowing grading by a collection of invertible objects instead of a single one. Secondly, higher comparison maps were defined by Sanders~\cite{Sanders13} in order to refine the analysis of the fibers of `lower' comparison maps, through an inductive process.
\end{Point}

\begin{Point}
The above comparison map is still very much concerned with the computation of~$\SpcK$. Moving away from this preoccupation, some first `geometric' results were established in~\cite{Balmer07}, like the decomposition of an object associated to a decomposition of its support, and applications to filtrations of~$\cat{K}$ by (co)dimension of support. These ideas naturally led to \emph{tensor-triangular Chow groups} in~\cite{Balmer13a} and further improvements by Klein~\cite{Klein16a} and Belmans-Klein~\cite{BelmansKlein17}, using the already mentioned relative tt-geometry of~\cite{Stevenson13}.
\end{Point}

\begin{Point}
In recent years, a great deal of progress followed from the development of the idea of \emph{separable extensions} of tt-categories. The ubiquity of this notion through stable homotopy theory, in connection with equivariant ideas, can be seen in~\cite{BalmerDellAmbrogioSanders15}. As a slogan, this theory extends tt-geometry from the Zariski setting to the \'etale setting. Implications for the spectrum are discussed in~\cite{Balmer16b}.
\end{Point}

\begin{Point}
Another area of tt-geometry which seems promising is the theory of \emph{homological residue field}, which aims at abstractly understanding the various `fields' which appear in examples: Morava $K$-theories, ordinary residue fields, $\pi$-points, etc. The reader can enter this ongoing project via~\cite{BalmerKrauseStevenson19,Balmer20a}.
\end{Point}

%%
%\begin{Point}
%We expect to write a more thorough treatise of abstract tt-geometry, beyond tt-classifications, in the hopefully not too distant future.
%\end{Point}

%----------------------------------------------------------------------
\newcommand{\etalchar}[1]{$^{#1}$}

%---------------------------------------------------------------------- 
%----------------------------------------------------------------------
\end{document}